\documentclass{amsart}
\usepackage{amssymb,amsfonts,amsthm}

\def\N{\mathbb{N}}

\def\Z{\mathbb{Z}}

\numberwithin{equation}{section}

\newcommand\mL{L\kern-0.08cm\char39}

\begin{document}
\title[Chaos in topological dynamics]
{Topological chaos: what may this mean?}
\author[F. Blanchard]{Fran\c cois Blanchard}
\address{Laboratoire d'Analyse et de Math\'ematiques Appliqu\'ees, UMR 8050 (CNRS-U. de Marne-la-Vall\'ee), UMLV, 5 boulevard Descartes, 77454 Marne-la-Vall\'ee Cedex 2, France}
\email{francois.blanchard@univ-mlv.fr}



\begin{abstract} We confront existing definitions of chaos with the state of the art in topological dynamics. The article does not propose any new definition of chaos but, starting from several topological properties that can be reasonably called chaotic, tries to sketch a theoretical view of chaos. Among the main ideas in this article are the distinction between overall chaos and partial chaos, and the fact that some dynamical properties may be considered more chaotic than others.
\end{abstract}

\maketitle

\section{Introduction}
In mathematics as elsewhere notions cannot go without names. Often those names are borrowed from ordinary language. They are picked up on account of some connection between their common meaning and the images raised by the mathematical definition. The association is mostly arbitrary: mathematicians often accept the first name that has been proposed because they are abstract-minded and ready to use any word as a conventional label for a property; their interest lies in the property, not in the name. Take, for example, `compactness': the mathematical notion slightly evokes the word. But nobody would say that a piece of Swiss cheese is compact, because there are holes in it, whereas the portion of the Euclidean space it occupies is compact mathematically, provided one considers it  closed, which does not offend common sense much. And what about manifolds, some of which barely evoke any of the meanings of the word?

For `chaos' the situation is not the same. The idea of chaos slowly emerged from experiments in physics. At some point scientists felt the need to express it mathematically. And then a misunderstanding appears. The physicist says `I think my intuition of chaos could be expressed by this mathematical property', never forgetting that this is, and must stay, questionable; the mathematician takes it for a convention and starts talking of chaos whenever the property is fulfilled. The physicist is expressing an opinion, the mathematician is freely using the word `chaos' as a label for a property or set of properties. 

In the meanwhile the latter forgets that there are several conflicting definitions of chaos. 
This creates some confusion in current mathematical literature: when one finds the word `chaos' in an article, according to context this may mean decrease of correlations, or positive Lyapunov exponents, or Li-Yorke chaos, or sensitivity to initial conditions, or Devaney chaos, or Hopf bifurcations. This is even more confusing for readers from other fields. In books the situation is not quite the same. There authors often acknowledge the fact that chaos cannot be summed up in one mathematical definition but should take into account several features, which shed different lights on the indeterminacy or disorder of the model that is being considered. 

What specialists mean by chaos is not the same in various sciences. Astrophysicists are often content to call `chaotic' a non-periodic dynamical behaviour. Unsatisfied with this elementary requirement theoretical physicists like Li, Yorke and Devaney proposed definitions which are commented below. In computer science people have started to think about chaos; they have not come up with definitions but there is a feeling that their intuition has to do with computational complexity, but that the idea of chaos ought to be distinguished from that of complexity.  Biology and system theory are likely to produce their own concepts of disorder. In their own fields mathematicians try to find properties that are evocative of what the word `chaos' may mean in other sciences, but  at no time can they give a definition, except by way of a convention devoid of any significance outside the field of mathematics. 

The aim of this article is limited. Focussing on topological dynamics, I would like to show --show in the ordinary, not the mathematical, sense-- that several of the notions from this theory make sense when talking about chaos and order, that they make more sense when not considered separately, and that, once this is accepted, one can distinguish degrees of chaoticity for dynamical systems. 

My purpose is not to introduce a new definition of chaos that would replace all others.  I want to say that chaos cannot be summed up in one formula, but there is another way of giving the word some sense in mathematics. In the past several people formulated tentative definitions of chaos. Each of them contains a part of truth, and they are not equivalent; as a consequence none of them should be taken for granted and none of them completely discarded. Those definitions are important, not because they are definitions but for two reasons: they give a form to intuitions; and, from an abstract point of view, they often indicate important correlations between factors of chaoticity inside particular classe of topological dynamical systems. In this article when saying that some of them are relatively weak I do not mean that they are no good. However weak, the Li-Yorke definition is still a kind of elementary requirement for chaos in any kind of dynamical system. The same is true about sensitivity in a different way. What should be emphasised is that 1. those definitions can be compared with each other and with some classical properties of topological dynamics that are usually regarded as chaotic without having ever been explicitly designated as `definitions of chaos'; 2. one may usefully distinguish various types of chaotic properties; 3. when talking about chaos one cannot by-pass the opposite notion, whether one calls it `determinism' or `order'. 

On the mathematical side the view of chaos I am tentatively developing here is certainly not comprehensive. The scope of the article is limited to topological definitions and techniques, which cover a rather minor part of all that has been written on the topic. Stochastic chaos is not taken into account, nor is the question of Lyapunov exponents; the reasons are that I am not especially competent in probability or differentiable dynamics, and that there is already enough to say about topological notions. Well inside the topological field, properties of attractors are not taken into account here; the explanation is in Section 2. Anyway there are plenty of holes in my documentation, some of which I am aware of and some of which I am not; the number of self-citations is due to the fact that I have been trying to understand something about chaos for many years, not to over-estimation of what I contributed. As the reader can see this paper, even though it refers to large sections of topological dynamics, is {\it not devised as a survey}. Its aim is to expound a coherent way of viewing chaos topologically. Whether the result is any good the readers will tell. What it is worth on the practical side they will tell too.

When an author claims that some mathematical property has to do with chaos what should be done about this notion, and what does he or she more or less expect people to do? A lot of things. Check which ones of the existing models possess this property and which ones do not; compare it with other chaotic properties; try to negate it and examine the outcome; consider how the definition may be strengthened or weakened; change quantifiers slightly in the definition and again examine the output. By the way some new properties are introduced, and can then be submitted to the same various manipulations. This is what I do here --partly. Suppose one does  this systematically, starting from five definitions: Li and Yorke's, Auslander and Yorke's, Devaney's, sensitivity and positive entropy. Then one is due to cover the whole field of topological dynamics --which some people occasionally call `the theory of chaos'-- introducing new properties, comparing them to previously known ones, hunting for examples, and by the way probably forgetting one's initial aim of representing chaos. To avoid this drift the scope of the article is limited, apart from some previous definitions of chaos, to a few notions that I believe to be essential in the theory as it stands, or empirically significant, or both; namely, equicontinuity, distality and weak mixing. 

Topological dynamics, like all other mathematical fields, developed along its own lines, which did not often  concern the meaning of the word chaos. For a more complete view of the theory the reader may consult books such as \cite{Br}, \cite{Au}, \cite{A}, \cite{dV}, \cite{G1}, \cite{Ku}; in addition many other notions that may be significant are described in a recent survey on local entropy \cite{GY}.  

Let me also mention an important point, which partly explains why the picture of mathematical chaos is so disorderly. It often occurs that {\it in a restricted class} of topological dynamical systems several chaotic properties or definitions of chaos coincide, or one implies another, while this is false when looking at the whole category of dynamical systems: among interval maps transitivity implies positive entropy, something which is trivially false elsewhere. It also occurs  that some of the properties which make sense in general never occur in a particular class. For instance cellular automata are never minimal, they are never transitive without being sensitive. Interval maps, symbolic systems and cellular automata have been particularly studied and many results about them are significant in terms of chaos. 

That each smaller class of dynamical systems has its own theory of chaos is important. The picture of chaos changes completely (while becoming more abstract) when one extends 
the scope from, say,  interval maps to the whole field of topological dynamics. Remark that it also changes everytime something new appears in one of the fields that are concerned. My own picture reflects the present state of affairs; step after step, new results and new examples, new models and new simulations will inevitably change the situation.

\smallskip
A few words are necessary about the peculiar character of this text. In all the article there is a mixture of mathematics and empirical ideas. In my mind this cannot be avoided when writing about chaos. 
As I said above
the notion of chaos does not come from mathematics. It was developed by physicists. Afterwards it gradually spread to computer science and other fields, where it took on slightly different meanings. Chaos can be interpreted in terms of mathematics, without ever becoming altogether a mathematical notion. For such an interpretation one should not forget that mathematical properties are relevant essentially when they are natural. What does this mean? In terms of physics, I would say that a property from topological dynamics is natural if there are dynamical systems with this property that may be used as models for physical phenomena some day. A property is natural for computer science if some systems with this property can be reasonably simulated by computations. This establishes a kind of hierarchy between dynamical properties that `feel chaotic' for the specialist: some of them are more natural, therefore more significant; sensitivity, positive entropy and their opposites appear to me as more natural than others. In the case of sensitivity, equicontinuity--an opposite of sensitivity--and positive entropy, there are systems with these properties that are used as models for experimental phenomena. This gives these properties, and their opposites, a particular importance empirically.  

The hierarchy of chaos is completely different from the one that arises when comparing the relative strength of properties in the mathematical field. This is illustrated by the fact that a property that is mathematically stronger than a chaotic property is not necessarily more chaotic than the latter. For instance completely scrambled systems (defined in Subsection \ref{LiYorke}) are not particularly more chaotic than Li-Yorke chaotic systems in general. 

The reader should never forget that all the stuff in the present article that is not downright mathematics is tentative and should be taken as such; and all the mathematical stuff, however rigorous, may be pointless for physicists or computer scientists. This does not put me in a position to obtain unambiguous results. Having two purposes that are orthogonal to each other, the article cannot be a complete success in the way a mathematical article may be when answering the questions that were at its origin. 
Here many properties are called chaotic, or deterministic, or neutral, because they evoke chaos or order in some way, others are called so for the sake of mathematical coherence, but none can be {\it proved} to be chaotic because there is no definition of chaos and there will not be any. The choices I made may be and should be discussed. The important point is that the word chaos has a meaning outside the field of mathematics. This meaning must be  reflected in the way mathematicians treat it. 

\medskip
The article is devised as follows. After the Introduction we first go through the five  `definitions of chaos', those of Li and Yorke, Auslander and Yorke, Devaney, sensitivity, and positive entropy; then we introduce a few other relevant properties from topological dynamics. There follows a short section about some intuitive features of chaos and order, and the methods that are used. In Section \ref{results} general results are reviewed. Section \ref{classes} is devoted to results that are valid in restricted classes and counter-examples in such classes. I finish with a section of provisional conclusions. I did not avoid repeating information  when it happened to be significant in several contexts.

\medskip \noindent {\bf Warning} 

Here when one reads `chaotic' this may have at least three different meanings. It may refer to the experimental acceptions of the word. At other times it means that a property has its place inside the topological theory of chaos I am trying to describe. Or else, that it has been considered chaotic by some author or authors. It is not always clear from context which one of the meanings is intended at each occurrence of the word. Here we are at the boundary of mathematics, which makes matters of vocabulary more complicated than they usually are.

The nouns `chaos', `disorder' and `indeterminacy' are used indifferently, and so are `order' and  `determinism'. Other authors give different meanings to these words, as can be seen in the title of  \cite{MS}. There `deterministic' means that the system is generated by some function (not by some random process), while the acception of `chaotic' is more or less the one I am using here. `Deterministic' is also the name of a precise property that was recently introduced in topological dynamics.  Contrarily, 
in ergodic theory a `deterministic' dynamical system is usually one with entropy 0. 

\section{Existing definitions of chaos and other properties}
Let me make it clear once more that in this section what is called a definition of chaos is one that has been explicity or implicitly coined as such by somebody. This does not mean that one must automatically accept it. The choice of further properties defined in this section is more arbitrary, but some explanations are given in due course; they may be so-called chaotic properties, or deterministic properties, or properties the character of which has to be examined. 

Here a few of those are introduced and discussed. 
For them to be valid technically it is sufficient to consider a topological dynamical system $(X,f)$,  and assume that $X$ is a metrisable space endowed with the metric $d$, and that $f$ is a continuous self-map of $X$; no assumption of compactness or surjectivity is strictly necessary. Nevertheless when one wants to explore the connections between various properties it is useful to assume that $X$ is compact.  That is what we shall do practically all the time. In all particular classes of dynamical systems that we consider in Section \ref{classes} the underlying space is compact, except the Besicovitch space for CA. In contrast  there are many non-surjective interval maps and cellular automata. 

So let $(X,f)$ be a dynamical system (for which $X$ is simply a metric space and $f\colon X\to X$ is continuous). By a {\it factor map} from $(X,f)$ to $(Y,g)$ one means a continuous map $\phi\colon X\to Y$ such that $\phi \circ f = g\circ \phi$, as usual in the theory; in this situation $(X,f)$ is called an {\it extension of} $(Y,g)$. 

\subsection{Definitions of chaos}
Let us first examine the five definitions of chaos  I mentioned above.

\bigskip
\noindent{\it Li-Yorke chaos.}

\smallskip
As far as I know the first topological definition of chaos ever introduced was Li and Yorke's. A subset $S$ of $X$
containing at least two points is called a \emph{scrambled set} if for any $x,y\in S$ with $x\neq
y$,

\begin{equation}\label{E:liminf}
\liminf_{n\to \infty} d(f^n(x), f^n(y)) = 0
\end{equation}
and
\begin{equation}\label{E:limsup}
\limsup_{n\to \infty} d(f^n(x), f^n(y)) > 0;
\end{equation}
in other words, the orbits of $x$ and $y$ get arbitrarily close to each other but infinitely many times they are at a distance greater than some positive $\epsilon$.
The system (or the map $f$) is called \emph{Li-Yorke chaotic} if it
has an uncountable scrambled set. 

The origin of this definition is in Li and Yorke's article \cite{LY}. Why should a system with this property be chaotic? That is not completely intuitive; one may discuss every single ingredient of the definition, but it turns out to make sense as it is. The first arguments come from the theory of interval transformations, in view of which it was most probably introduced. For such maps the
existence of one scrambled pair implies the existence of an uncountable
scrambled set \cite{KuS}, and it is not very far from implying all other properties that have been called chaotic in this context; for a survey of known results see \cite{Ru}. As we shall see in the next section, in general the Li-Yorke property passes several theoretical tests for being an elementary requirement for chaos: Li-Yorke chaos has been proved to be a necessary condition for many other  `chaotic' properties to hold.

\bigskip
\noindent{\it Sensitivity.}

\smallskip
Sensitive dependence on initial conditions (sensitivity for short) is rather intuitively a chaotic property. Nowadays it is usually defined as follows: $(X,f)$ is said to be {\it sensitive} if there exists a sensitivity constant $\eta>0$ such that for any $x \in X$, any $\epsilon>0$, one can find $y$ with $d(x,y) <\epsilon$ and $n$ such that $d(f^n(x), f^n(y))\ge \eta$. The initial idea is due to Lorentz.  Originally in \cite{Gu} Guckenheimer  introduced this definition for interval maps, with a further requirement that the property hold  on a set with positive Lebesgue measure. Later authors preferred the stronger definition above, which extends to all metric spaces and does not contain any probabilistic assumption. 
Like Li-Yorke chaos, sensitivity is a consequence of several other chaotic properties, but it is not equivalent by far. 

Heuristically the idea is a very interesting one. Here is an elementary relation between metric and predictability. In a sensitive system with sensitivity constant $\eta$, however close a point $x$ is to a point $y$, even if one can compute the orbit of $y$ precisely it is impossible to predict with precision better than $\eta$ how the orbit of $x$ behaves in the long run under iterations of the map. Even among non-scientists this property is widely known as  `the butterfly effect', certainly in part because of the poetical flavour of the phrase. 

\bigskip
\noindent{\it Auslander-Yorke chaos.} 

\smallskip
Later Auslander and Yorke \cite{AY} borrowed the strong, purely topological, version of sensitivity and introduced their definition of chaos by associating it with transitivity. This is called {\it Auslander-Yorke chaos}. By {\it transitivity of} $(X,f)$ one means that for any two non-empty open subsets $U$ and $V$ of $X$ one can find an integer $n>0$ such that $U\cap f^{-n}V \ne \emptyset$. A transitive map is necessarily onto; when $X$ is compact transitivity is equivalent to the existence of a dense orbit, that is, that of a point $x \in X$ such that the set $\{x,f(x),\hdots,f^n(x),\hdots\}$ is dense in $X$. 

The Auslander-Yorke definition combines one property, sensitivity, which is strongly evocative of indeterminacy or disorder, and another--transitivity--which is not; the second point is explained in Subsection \ref{further}. However, there are arguments for assuming transitivity when trying to understand chaos: usually it is a simplifying hypothesis; one may choose to look first for notions of chaos that are valid for systems having a certain amount of homogeneity under iterations of $f$; finally, unless the contrary is established one may assume that any feature of chaos should appear among transitive systems.  

\bigskip
\noindent{\it Devaney chaos.}

\smallskip
Still later Devaney introduced his own definition of chaos \cite{D}; in the sequel it is called {\it Devaney chaos}. It relies on three assumptions. A map is called chaotic in his sense if it is sensitive, transitive and if $X$ contains a dense set of $f$-periodic points: in other words, if it is Auslander-Yorke chaotic with  dense periodic orbits. Not much later Banks et al. \cite{Ba} remarked that transitivity and dense periodic orbits imply sensitivity. The same was proved again, assuming compactness, by Glasner and Weiss \cite{GW}; apart from this result they did not know to have been already proved, their article contains several interesting remarks about chaos and results. 

My interrogations with regard to Devaney chaos are explained at length in Subsection \ref{AYD}.

\bigskip
\noindent{\it Positive entropy.}

\smallskip
In \cite{F}, without using the word chaos Furstenberg chose to call `deterministic' all compact dynamical systems with topological entropy $0$. Later Glasner and Weiss \cite{GW}, in a discussion of Devaney's definition, proposed positive entropy as the essential criterium of chaos.

Positive topological entropy is a strong property. Again consider a topological dynamical system $(X,f)$ with $X$ compact. For any finite cover
${\mathcal C}$ of $X$, let ${\mathcal N}({\mathcal C})$ be the
minimal cardinality of a sub-cover of ${\mathcal C}$. Let ${\mathcal C}^n={\mathcal C}\bigvee f^{-1}{\mathcal C}\bigvee ... \bigvee
f^{-(n-1)}{\mathcal C}$, and
$c_f({\mathcal C},n) = {\mathcal N}({\mathcal C}^n)$.
The exponential growth rate of the function $c_f({\mathcal C},n)$ is the topological entropy of
${\mathcal C}$, and the {\it
topological entropy\/} $h(X,f)$ of $(X,f)$ is the supremum of the
topological entropies of finite open covers. Bowen introduced other definitions in the non-compact setting \cite{Bo}; they are equivalent to the one above in the compact case.

Regarding positive topological entropy as a chaotic feature finds a heuristic justification in ergodic theory, or more precisely in the part entropy plays in information theory. There measure-theoretic entropy is interpreted as a measure of the indeterminacy of an invariant process: a process with discrete states is called deterministic if its past states almost surely determine its state at time 0, which is equivalent to its entropy being null. The link with topological entropy is the Variational Principle: the topological entropy of the system $(X,T)$ is the supremum of the entropies of all $f$-invariant probability measures on $X$ (for more on this topic see \cite{Wa}). Therefore a deterministic topological system in the sense of Furstenberg is one such that any $f$-invariant probability measure on $X$ has entropy 0.  

\subsection{Further properties and a first discussion}\label{further}

\smallskip
A really convincing theory of topological dynamics exists only with the assumption that the space $X$, in addition to being metric, is also {\bf compact}. Most general results concerning chaotic properties, but not all, are obtained under this hypothesis. 

 Intuitively {\bf surjectivity} evokes the physical notion of equilibrium (but the latter is richer).
Surjectivity does not look an essential assumption when dealing with chaos. Keep in mind the fact that transitivity implies surjectivity. To the contrary Li-Yorke chaos, sensitivity, entropy make sense even when the transformation is not onto. In topological dynamics it is often expedient to assume that the action is surjective, but one has always the option to restrict the action of $f$ to the $\omega$-limit set of $(X,f)$, by the way losing some information about the overall behaviour of the transformation. Still, when $f$ is not onto one has to make sure that Li-Yorke chaos or sensitivity also occur in the $\omega$-limit set, which is not always the case. 

It would be undue to call {\bf transitivity} a chaotic property. This is best illustrated by irrational rotations of the circle: they are transitive maps, while being isometries, which prevents them from satisfying any of the popular criteria of chaoticity; moreover they are adequate models for the movement of a bicycle wheel around its centre, which is among the most orderly phenomena one may observe. Anyway, transitivity together with dense periodic orbits imply sensitivity, and that is a chaotic property even if a very weak one. Note that if transitivity is hardly connected with chaos, non-transitivity intuitively involves some degree of order: starting from some neighbourhood there are regions of the phase space that one can never reach. Thus transitivity, without being a really chaotic property, prevents the occurrence of some elements of rigidity among dynamical systems. There is a short survey of transitivity in relation with chaos in \cite{KoS}.

Minimality is a stronger property. A dynamical system $(X,f)$ is called  {\it minimal} if there exists no closed $f$-invariant subset $S \subset X$ except $X$ and  the empty set. This property makes sense mostly in the compact setting: any point of a minimal compact system has a dense orbit. As the reader may check throughout the article, minimality appears to be transverse to chaos. Its main interest for us is that it forbids the existence of any periodic orbit when $X$ has no isolated points; it also plays an essential part in the disjointness of two systems, defined below in Section \ref{methods}.

Abstractly, having {\bf dense periodic orbits} is not much more of a chaotic feature than transitivity; this is illustrated in symbolic dynamics (Subsection \ref{symb}) in \cite{W} and even more strikingly in \cite{BH}, and also in the field of  cellular automata, where both the most and the least chaotic maps have been proved to have dense periodic points (Subsection \ref{CA}). Nevertheless I cannot think of any ordered physical phenomenon that could be modelled by a system having dense periodic orbits, and this makes me believe that Devaney's intuition has some practical value. Having dense periodic points looks a chaotic property, provided it goes along with some other reasonable assumption.
Judging that dense periodic orbits are too specific a requirement, Glasner and Weiss suggest to replace it by a weaker one, existence of a dense set of minimal points, or by that of an $f$-invariant measure with support, which is even weaker and, combined with transitivity, implies that the system is either minimal equicontinuous or sensitive  \cite[Theorem 1.3]{GW} ; for them those are rather technical assumptions, while they consider chaos as coinciding with, or defined by, the existence of a positive-entropy measure. 

Here we must examine one classical property that has never been chosen as a definition of chaos: {\bf weak mixing}. Weak mixing deserves to be introduced in this paper, not because there are many weakly mixing models of phenomena --I do not know any that has not got positive entropy at the same time-- but on account of its special position in the theory:  it implies both Li-Yorke chaos \cite{I} and, trivially, sensitivity;  it is much stronger than transitivity, and weaker than positive entropy in some way (this is explained in Subsection \ref{LiYorke}); finally it plays a r\^ole in some disjointness theorems also involving deterministic properties. It is defined as transitivity of the Cartesian product system $(X\times X, f\times f)$; long ago Furstenberg discovered that this property is equivalent to transitivity of all finite Cartesian powers $(X\times X \times \ldots \times X, f\times f \times \ldots \times  f)$ \cite{F}: so it is not merely a property of pairs of neighbourhoods, but jointly concerns any number of neighbourhoods. This illustrates the fact that weak mixing is really much stronger than transitivity. Again most of its significance and consequences disappear when one drops the compactness assumption on $X$. Then it still implies sensitivity, but I doubt that any of its other consequences remains true.

\medskip
Let us finally define two strong `deterministic' properties. The first, equicontinuity, contradicts sensitivity and Li-Yorke chaos, the other one, distality, contradicts Li-Yorke chaos only, and both play an important part in general topological dynamics. 

An {\it equicontinuity point} $x$ of $(X,f)$ is one such that for any $\epsilon>0$ there is $\eta>0$ such that if $y \in X$ with $d(x,y) \le \eta$ then $d(f^n(x),f^n(y) \le \epsilon$ for any power $n\ge 0$. The system $(X,f)$ is called {\bf equicontinuous} if for any $x \in X$, any $\epsilon >0$ there is $\eta>0$ such that $d(x,y) <\eta$ implies that $d(f^n(x), f^n(y)) <\epsilon$ for any integer $n$. On a compact space a system is equicontinuous if and only if all its points are equicontinuity points. There are systems with equicontinuity points, even transitive ones, that are not equicontinuous \cite{GW,AAB,Ku}.
$(X,T)$ is called {\bf distal} if whenever $x \ne y$, there is $\delta>0$ such that $d(f^n(x),f^n(y) \ge \delta$ for any $n$. Equicontinuity implies distality; this is not completely obvious. 

Both properties have been studied extensively. Their importance in the theory is developed in \cite{Au,dV}. Equicontinuity forbids any kind of sensitivity. More precisely in an equicontinuous system  there exist no sensitivity points, that is, no points $x$ such that there exists $\epsilon>0$ such that for $\eta>0$ one can find $y\in X$ and $n$ with $d(x,y)<\eta$ and $d(f^n(x), f^n(y))\ge \epsilon$. Intuitively equicontinuity is a deterministic property, and so is its weaker version, having equicontinuity points. The antagonistic relation between various equicontinuity properties (among them, existence of equicontinuity points) and sensitivity among transitive maps is described in \cite{AAB}. Distality can be considered deterministic because it implies entropy 0, actually forbidding the existence of any non-empty scrambled set, and also on account of one of Furstenberg's disjointness theorems:  weakly mixing systems are disjoint from all minimal distal systems \cite{F}.

\subsection{Properties that are not covered.} 

For different reasons some chaotic or deterministic properties are not examined in this paper.

\smallskip
One might, probably should, consider many other topological properties, like strong or mild mixing, scattering, rigidity etc. but in this paper my purpose is to concentrate on a few especially significant ones. Presently none of the existing examples of strongly or mildly mixing systems with entropy 0 is likely to be used for representing an experimental phenomenon, while the two properties imply weak mixing.  So, better consider only one of those abstractly chaotic properties, here weak mixing, and leave others for future investigation. For similar reasons I decided not to insist on almost distal systems; while excluding   Li-Yorke chaos, at present almost distality has not got the same importance as distality in the theory.  

Still another matter is that of attractors and $\omega$-limit sets. When $X$ is compact the $\omega$-limit set of $(X,f)$ is just the intersection $\bigcap_{n=1}^\infty f^n(X)$, that is, the maximal attractor.  A small limit set, or the existence of one or several small attractors inside $X$, can be regarded as features of order. The topological properties of $\omega$-limit sets are  important in some particular settings, for example in one-dimensional dynamics where the fact that the system is attracted by a Cantor is correlated with chaotic properties. 

I chose not to treat properties of attractors for two reasons. The first objection is that in various other classes of dynamical systems one can construct systems $(X,f)$ having a complicated $\omega$-limit set on which the transformation acts as the identity: the system may exhibit some unpredictability features in the way it evolves {\it towards} equilibrium, but certainly not in the way it behaves {\it at} equilibrium. This happens among cellular automata and also in the differentiable setting. For the latter situation, that of so-called non-chaotic strange attractors, see for instance \cite{GJK}, in which Lyapunov exponents are used. 
 Anyway, this shows that the form of an attractor may have no meaning whatsoever with regard to chaoticity of the map acting on this attractor. My second reason is that taking attractors into account would double the size of the article, and the matter could not be treated suitably without considering differentiability and Lyapunov exponents. The disadvantage is that one leaves aside any feature of chaoticity that might occur during transition to equilibrium. Of course I would not say that attractors have to be discarded when dealing with chaos. 

As I already said the probabilistic side is not addressed here. Of course it should be taken into account if one wanted to review the many mathematical attempts to describe chaos, and there are basic scientific reasons for addressing this field when trying to understand disorder. But here my choice is to limit ourselves to topology in a broad sense. It is for the same reason that I do not treat Lyapunov exponents. In differentiable dynamics the existence of a positive Lyapunov exponent is  considered a chaotic feature. The relation with positive entropy is certainly a good argument \cite{Y}. But one cannot define the Lyapunov exponents when there is no geometric structure. 

\section{Postulates and methods}\label{methods}

Let me first introduce a distinction that I believe to be essential. When comparing
 attentively topological properties that have been called 
 or are often considered chaotic, one remarks that there are two kinds. Some, which I call {\bf overall (chaotic) properties}, hold everywhere in $X$; among them, sensitivity (in any neighbourhood of any point $x$ one can find $y$ such that, etc.), Devaney chaos (because both transitivity and dense periodic orbits are overall properties), weak mixing, strong mixing and more. The systems they distinguish can be said to be everywhere chaotic, in a sense which can be weaker or stronger according to the property. On the other hand, definitions like those of Li-Yorke chaos and positive entropy say that a system contains some kind of chaos, but not necessarily everywhere; let us call them {\bf partial (chaotic) properties}. 
To see whether some property belongs to the first or second family one has to consider occurrences of the quantifiers $\forall$ and $\exists$ in its definition. But this must be done cautiously. The point is not whether the definition begins with $\forall$ or $\exists$, as I wrote unthinkingly somewhere. The definition of sensitivity as given above begins with `there exists $\epsilon>0$\ldots', but afterwards the requirement appears to be for all points, and this makes it an overall property. 

Of course some properties of dynamical systems enter with difficulty, if at all, into one of those two classes. One of those that are mentioned in this article is the existence of a second-category scrambled set. 

In topological dynamics properties are much easier to manipulate when they are preserved either under factor maps or under lifts, or even better under both. A property which is neither is a rather loose one and not very significant. In the present context overall chaotic properties ought to be preserved under factor maps if one wants them to tell something important about the structure of the system. This is true for all those that are mentioned in this article except for sensitivity. But properties that are significantly chaotic are not necessarily nice tools in the theory. On the other hand, partial chaotic properties cannot be expected to be stable under factor maps. If a system contains some amount of chaos any larger system, or any system that contains it structurally, should also have some chaoticity of the same kind, both for commodity in the theory and for the intuitive reason that partial chaos should survive embedding (if not shrinking). In other words, partial properties should be preserved under unions, which they usually are, and under extensions, which is often harder to prove: see some of the questions at the end of the article. Note that Li-Yorke chaos is not stable under extensions (Example 1), while positive entropy is. 

There are correspondences between overall properties and partial properties. Given one overall property it is always possible to find a corresponding existence property, often several. The reverse is true too (but in the case of sensitivity a property like `sensitivity somewhere' is so weak that it has very little any significance). Once this is done there remains to investigate all options and retain the most significant one -or ones.  

\medskip
Another important point is that there are antagonistic relations between chaotic properties, whether overall or not, and deterministic or order properties. Given one chaotic property, one always finds some opposite property that makes it completely impossible for the first to hold.  Most of the time one chaotic property can be contradicted in several different ways, which means that the choice is not completely obvious. Anyway if one wishes to really understand the meaning of some possibly chaotic property it is necessary to formulate at least one opposite definition.   

Just like for chaotic properties one may consider overall and partial order properties. Since we wish the order property to contradict the chaotic property, when starting from a partial chaotic property one obtains overall deterministic properties, and the other way round when the chaotic property is overall. Good examples of partial deterministic properties are the existence of equicontinuity points, or that of a non-trivial equicontinuous factor, in a system. They may co-exist with partial chaos in the same dynamical system. For example there are systems with positive entropy and a non-trivial equicontinuous factor, that is, chaotic systems having a heavily deterministic element of structure. By the way, as is the case for overall chaotic properties, overall deterministic properties are more significant when they are preserved under factor maps.

I do not examine partial order properties very consistently here, essentially because up to now they have not been explored very much in topological dynamics. What would be a partial property corresponding to entropy 0? There are different sensible answers, and I cannot tell which one is the most significant. In the motivation of this article there is a bias in favour of chaos. Maybe future research will compensate this bias, or else somebody will find good arguments for an absence of symmetry between chaos and order. 

\medskip
Let us now consider a theoretic tool : disjointness. Two dynamical systems $(X,f)$ and $(Y,g)$ are said to be {\it disjoint} if the unique closed, $f\times g$-invariant subset of $X\times Y$ that projects to $X$ and $Y$ respectively is $X\times Y$ itself. For two systems to be disjoint, one of them at least must be minimal. This is a limitation: for instance it is true that all weakly mixing systems are disjoint from all {\it minimal} distal systems, but it is false that all weakly mixing systems are disjoint from all distal systems. Of course it is hopeless to prove a disjointness result between a partial chaotic property and some deterministic one; both properties involved must be overall. 

Disjointness is a very abstract property. One may consider two disjoint systems as having very different dynamical structures, and two disjoint properties, like weak mixing and minimal distality in the claim above, as being contradictory in a rather deep sense. This view is partly misleading: there is evidence in the theory that two disjoint dynamical systems may have deep features in common. We can nevertheless take disjointness as a weak test of the non-existence of common features for two topological dynamical systems; right now I do not know of any stronger one in this setting. Disjointness may also be considered as a test of significance for both properties involved: I am using it that way in this article. There are properties that are not disjoint from any other property, as is illustrated in the next section by completely scrambled systems. These properties can never be very strong. 

\medskip
Finally, the reader should be warned that chaos and order are not clear-cut notions. As can be seen below, a non-transitive system has some features of order, whereas one can hardly consider transitivity to be a chaotic property; a system that is not Li-Yorke chaotic is a very ordered one, but one cannot consider Li-Yorke chaos as involving much chaoticity.

\section{General results}\label{results}
At some point we must really get into technicalities. That is what we do now. In this section I consider first Li-Yorke chaos, then sensitivity, Auslander-Yorke and Devaney chaos, weak mixing and finally positive entropy, in relation with other topological properties, whether chaotic or deterministic, and in the end I sketch some scales of chaos. 

Almost all the results stated in this section are valid only when $X$ is compact. Then let us take compactness as a standing hypothesis.

\subsection{Li-Yorke chaos}\label{LiYorke}

Here I am applying the prescriptions listed in the section above: comparing the Li-Yorke definition with other definitions existing in the field, strengthening and weakening it, finding opposite properties, and doing the same all over again with properties newly introduced; but not too systematically, in order to keep the article within reasonable limits and refrain from boring the readers.

The Li-Yorke condition is obviously weak. One cannot expect scrambled sets to tell all about the dynamics of a transformation: their definition relies on pairs being scrambled or not, that is, on the behaviour of joint orbits, independently of what occurs in the neighbourhood. Stronger properties like weak mixing or positive entropy can be expressed in terms of the behaviour of neighbourhoods of pairs, and certainly properties of joint orbits cannot be sufficient to characterise them. 

Li-Yorke chaos is a partial property. According to the list of requirements in the previous section it would be better if any extension of a Li-Yorke chaotic system were also Li-Yorke chaotic; actually this is false: 

\medskip \noindent {\bf Example 1} (Blanchard and Huang): the article \cite{BDM} describes a minimal symbolic system $(Z,\sigma)$ which is a bounded-to-one extension of an adding machine; it contains only bounded scrambled sets, and in particular a recurrent scrambled pair $(z,z')$. 

Let $(Y,\sigma\times \sigma) \subset (Z\times Z,\sigma\times \sigma) $ be the orbit closure of $(z,z')$ under $\sigma\times \sigma$: this is easily seen to be a transitive system having only bounded scrambled sets, and it is not minimal since it contains the diagonal $\Delta_Z$. Define a factor map $\pi\colon (Y,\sigma\times \sigma) \to (X,f)$ by collapsing all points of $\Delta_Z$ to one fixed point. $(X,f)$ is transitive as a factor of $(Y,\sigma\times \sigma)$, and contains a periodic point: this implies Li-Yorke chaos by  \cite{HY2}. Thus $(Y,\sigma\times \sigma)$ is a non-Li-Yorke chaotic  extension of the Li-Yorke chaotic system $(X,f)$.

\medskip
Let us compare Li-Yorke chaos with other definitions. It is proved in \cite{HY2} that Devaney chaos implies Li-Yorke chaos--actually that Li-Yorke chaos results from transitivity and the existence of just one periodic orbit. The existence of a periodic orbit is crucial, since transformations like the shift action on Sturmian systems or the Morse system are minimal, hence transitive, and sensitive. So Auslander-Yorke chaos does not imply Li-Yorke chaos. In the other direction Li-Yorke chaos, which is not an overall definition, cannot imply the overall property of sensitivity: think of  the union of a Li-Yorke chaotic system and an equicontinuous one. A more subtle, transitive counter-example results from \cite{AAB}:

\medskip \noindent {\bf Example 2:} Theorem 4.2 of \cite{AAB} states the existence of a transitive system $(X,f)$ having equicontinuity points,  which by Theorem 4.1 of the same article contains a fixed point; the existence of equicontinuity points in $(X,f)$ prevents sensitivity, while by \cite{HY2} transitivity and a fixed point imply Li-Yorke chaos.  

\medskip 
Another test of the value of the Li-Yorke assumption is that all systems for which it is known {\it not to hold} have strong deterministic features: to quote only the most conspicuous, their entropy is 0 and they cannot be weakly mixing. Whenever there is no uncountable scrambled set few other chaotic properties may hold. 

In order to strengthen Li-Yorke chaos one may ask the scrambled set $S \subset X$ to be not merely uncountable but , say, big topologically inside $X$: second category, or residual, or equal to $X$ itself. In \cite{BHS} the authors try to gather everything that is known about the size of scrambled sets, and solve several related open questions. Some results on the size of scrambled sets and some examples of second-category or residual scrambled sets are given: let me mention that a minimal system never has a residual scrambled set; in the other direction, a weakly mixing system always has a second category scrambled set if one assumes the Continuum Hypothesis. 

When $X$ is a scrambled set the system is called {\it competely scrambled} in \cite{HY1}, where some properties of those systems are given--in particular such a system is never minimal--and one example is described; in \cite{HY3} the same authors prove the existence of a completely scrambled system which is also weakly mixing. Being completely scrambled is an overall property corresponding to Li-Yorke chaos, and at first sight a strong one. In fact it is not so strong, and not a very chaotic property. Again in \cite{BHS} it is shown that a completely scrambled system may have a non-trivial factor without any scrambled pair. Without any further information one would expect completely scrambled systems to be disjoint from minimal distal systems. This is false: given any minimal distal system, \emph{a fortiori} given any minimal equicontinuous system, one constructs a non-disjoint completely scrambled system \cite{BHS}. 

Distributional chaos was introduced in \cite{SS} as another way of strengthening Li-Yorke chaos into a more significant property. Other non-equivalent versions were proposed later. This was the origin of wide investigations, reviewed in \cite{S}, but there are still many points that should be addressed.

Let us point out that an extreme opposite of Li-Yorke chaos is certainly the absence of scrambled pairs. Equicontinuous systems and the larger class of distal systems have no scrambled pairs, hence the significance of the aforementioned result on non-disjointness.  There are further examples in \cite{B-M}. Systems without scrambled pairs have been explored in \cite{B-M}, where they are called {\it almost distal}. In the same article it is shown that a factor of an almost distal system is almost distal, which implies that any extension of a system having scrambled pairs also contains scrambled pairs. Weaker properties that still preserve a strong deterministic flavour were introduced by various authors, they are reviewed in \cite[Section 4]{BHS}. 

Weak mixing and one corresponding partial property, partial weak mixing, are dealt with in Subsection \ref{wm}; positive entropy and one corresponding overall property, u.p.e., are examined in Subsection \ref{posent}. Once it is admitted that weak mixing and partial weak mixing, u.p.e. and positive entropy are chaotic properties, one can re-examine how the size of scrambled sets varies according to those properties. A very tentative observation is that the scrambled sets of systems with the strongest chaotic properties are big but not too big. Weak mixing, and therefore u.p.e., imply the existence of a second category scrambled set, if one accepts the Continuum Hypothesis \cite{BHS}. The same article gathers examples of systems with positive entropy having at most first category, or at most second category, or at most residual scrambled sets. But the full shift and many related weakly mixing positive-entropy subshifts have no residual scrambled sets, and the same is true for all minimal positive-entropy systems. Contrarily completely scrambled systems have been proved to have entropy zero \cite{BHR}, which again suggests that the property is not a strongly chaotic one. 

My conclusion is that Li-Yorke chaos prevents a system from being essentially ordered, but does not guarantee that it is more than superficially chaotic. 

\subsection{Sensitivity}

At present there is less to say about sensitivity. It is one of the assumptions of Auslander-Yorke and Devaney chaos, and we already know that it cannot be compared with Li-Yorke chaos. It is an overall property, and it is not hard to formulate significant corresponding partial properties: existence of points that are not equicontinuous, or of a $G_\delta$ of sensitive points for some sensitivity constant, etc; such properties are extremely weak. Though an overall property, sensitivity is not stable under factor maps. 

An obvious opposite of sensitivity is equicontinuity, as was remarked above; furthermore a transitive system that is not sensitive has equicontinuity points \cite{GW}. Distality is not an opposite assumption, since any minimal distal system which is not equicontinuous is sensitive \cite{AAB}. Weak mixing implies sensitivity: this is an immediate consequence of the density of scrambled pairs in the Cartesian square of a weakly mixing system. 

I claimed that sensitivity is weak. For instance any 0-entropy system $(X,f)$  has a 0-entropy symbolic extension, which is often sensitive while obviously strongly related to $(X,f)$ \cite{BD}. Sturmian systems (Subsection \ref{symb}) and Toeplitz extensions of adding machines are classical examples. On the other hand some property resembling sensitivity but stronger, which for instance would not be satisfied by very un-chaotic systems like Sturmian shifts or the Morse shift, would be a better test of chaos. Finding one such would be an achievement; presently I do not know how this could be done. 

There are also interesting developments about sensitivity and weak mixing in \cite{YZ}. 

\subsection{The Auslander-Yorke and Devaney definitions}\label{AYD}

In general, when combined with transitivity, sensitivity does not get much more chaotic than it is by itself. One sees in the next section that in the class of subshifts, which exhibit a large range of properties,  Auslander-Yorke chaos does not mean much. In classes with a smaller range of properties, like interval maps and CA, it can be considered as a significant requirement for chaos. 

What may be the opposite of Auslander-Yorke chaos? Looking for an overall deterministic property, I would say equicontinuity and transitivity. Equicontinuity alone would also do. Remark that  Auslander-Yorke chaotic systems and minimal equicontinuous systems are not always disjoint, as shown by Sturmian systems (Subsection \ref{symb}). 

Devaney's definition is certainly more significant, since it implies Li-Yorke chaos, but it describes a very special kind of chaos. 

It is easy to find an opposite of Devaney chaos. Again it is an overall property, so that equicontinuity and transitivity would do since together these properties imply minimality, hence no periodic orbits. The paradox is that this is exactly the same opposite as for Auslander-Yorke chaos. Again Devaney chaotic and minimal equicontinuous systems are not always disjoint; there is an example in \cite{BH}.

Let us anticipate on Section  \ref{classes}  and consider briefly what happens among interval maps, cellular automata and subshifts. 
In the setting of interval maps a dense set of periodic points is strongly correlated with all chaotic phenomena; the many results substantiating this claim are reviewed in  \cite{Ru}. Cellular automata are dynamical systems that are easy to simulate computationnally. They divide into two classes, those having equicontinuous points and those that are sensitive. All left- or right-closing CA have dense periodic orbits, and these two classes contain essentially cellular automata having broadly chaotic properties. On the other hand all surjective CA with equicontinuous points also have a dense set of periodic points. Thus dense periodic orbits cannot be regarded as a chaotic feature among CA. 

In symbolic dynamics, another field having strong links with computer science, the picture is different and more complicated. Suppose we restrict ourselves to subshifts having a simple combinatorial description, which can thus serve as easy models for physical phenomena: then transitive systems that are chaotic on many accounts have dense periodic orbits, and those that are rather non-chaotic have none at all. But when considering the whole field one finds subshifts that are very chaotic without having dense periodic orbits, and systems that are not very chaotic but with dense periodic orbits. These coincidences of properties may look pathological to somebody who is familiar with interval maps only. 

However frequently they appear in actual models, in the general theory dense periodic points are an out-of-the-way requirement, and not a chaotic one in some settings. 
Topological dynamists will not admit that a minimal system with positive entropy (which does not contain any periodic orbit) is not chaotic in some sense; or that a 0-entropy system with dense periodic orbits is chaotic, if they consider positive entropy as the sole real test of chaos. On the other hand, the various existing examples of minimal positive-entropy systems are not very easy to construct and certainly harder to use as models for physical phenomena. Here the view of mathematicians
may diverge from that of other scientists. 

\subsection{Weak mixing}\label{wm}

Weak mixing is an overall property. It implies Li-Yorke chaos and actually more. A weakly mixing system always contains a dense uncountable scrambled set \cite{I}: to my knowledge this result is the first one that was proved about Li-Yorke chaos in the general setting. 

Weak mixing strongly contradicts distality since a weakly mixing system has a dense uncountable scrambled set, while a distal system has no scrambled pairs at all. There is also the already quoted  Furstenberg theorem proving disjointness between weak mixing and minimal distal systems. 

It would be nice and useful to display some partial chaotic property corresponding to weak mixing. In \cite{BH} the notion of weakly mixing set, derived from a result in \cite{XY}, is introduced. $A$ is a
weakly mixing subset of $X$ if and only if for any $k\in \mathbb{N}$, any choice of non-empty
open subsets $V_1,\cdots, V_k$ of $A$ and 
$U_1,\cdots,U_k$ of $X$ with $A\cap U_i\neq \emptyset,
i=1,2,\cdots,k$, there exists $m\in \mathbb{N}$ such that $T^{m}V_i\cap
U_i \neq \emptyset$ for every $1\le i \le k$.
This definition is rather complicated and still being studied; simpler versions may be found some day. A system $(X,f)$ is called {\it partially weakly mixing} if $X$ contains a weakly mixing subset; unsurprisingly, weak mixing is equivalent to the fact that the whole set $X$ is a weakly mixing set. It is also proved that positive entropy implies partial weak mixing, which in its turn implies Li-Yorke chaos. 

Partial weak mixing is a good candidate for being the partial notion corresponding to weak mixing. It is still unknown whether it is preserved under extensions. But the notion of sensitive sets defined in \cite{YZ} is another candidate. Both notions were introduced quite recently and their relations have not yet been investigated. 

\subsection{Positive entropy}\label{posent}

A system with positive entropy is usually considered chaotic by mathematicians. Topological dynamical systems with entropy 0 were called deterministic in \cite{F}; the reason is that in ergodic theory measure-theoretic entropy is usually interpreted as a measure of indeterminacy, and entropy 0 as a test of determinism; finally the Variational Principle allows one to transfer those notions to topological dynamics. 

In \cite{B-M} the authors proved that positive entropy implies Li-Yorke chaos. 
Positive entropy relies on the existence of a finite cover, the complexity of which grows exponentially under iterations of $f$; this definition makes it a typical partial property. For some dynamical systems any finite cover by non-dense open sets has positive entropy. This property, which is called {\it uniform positive entropy}, u.p.e. for short, is an overall property corresponding to positive entropy. It is a strong one, since any system with u.p.e.  is disjoint from any minimal 0-entropy system \cite{Bl}. It is preserved under factor maps. Moreover u.p.e. implies weak mixing.
Another possible overall property is completely positive entropy, in other words, no non-trivial zero-entropy factor exists. It is weaker than u.p.e., and though it is preserved under factor maps, a rather loose property. 

\subsection{Degrees of chaos}\label{degrees}
There are so many relations between the chaotic properties listed above (perhaps putting aside Auslander-Yorke chaos) and their opposites that basing an abstract theory of chaos on them looks reasonable. I do not claim that it is relevant. 

It certainly occurred previously to many scientists that among chaotic models some look more chaotic than others. In view of the results listed above, in particular those concerning implications or disjointness, it is possible to draw different scales of chaoticity and determinism. Such a scale is a short list of properties, each of which implies the next. But in order for this to be really meaningful, ideally one should have packs of three scales instead of an isolated one: one scale of partially chaotic properties, one scale of corresponding overall chaotic properties, and an inverse list of order properties that are opposites of the chaotic ones. When this is possible disjointness theorems involving one overall chaotic property and its opposite would testify that the chaotic property is not too soft. 

Ideally there should be only one way to do this, and the three scales should satisfy all my requirements. Actually there are several possible threefold scales and none of them is completely satisfactory. From a purely mathematical point of view this is a shame. If one considers this as a manifestation of the reality principle it is rather reassuring. 

As explained above sensitivity and Li-Yorke chaos are two elementary, non-equivalent requirements for chaos; formally sensitivity, as an extremely weak overall property, is slightly less satisfactory but it is quite significant empirically. 
Since sensitivity and Li-Yorke chaos cannot be compared, such scales may contain at most one of these two properties. None of them has a deep significance. They are here to represent elementary requirements for chaoticity in comparison with stronger properties. 

With Li-Yorke chaos as the lowest degree, one may draw two parallel scales, the first concerning partial chaos:
$$\hbox{positive entropy } \Rightarrow \hbox{ partial weak mixing }\Rightarrow \hbox{ Li-Yorke chaos}$$ and 
the other one containing overall notions only, each of them implying the corresponding one in the first scale: $$\hbox{uniform positive entropy } \Rightarrow \hbox{ weak mixing }\Rightarrow \hbox{ dense  scrambled set;}$$
the existence of a dense scrambled set is a very weak overall property, but one cannot replace it by the existence of a residual scrambled set, which is not implied by weak mixing and may not be a really chaotic property. Investigations on what property one might put in its place would be welcome. 

A tentative inverse scale of determinism would be
$$\hbox{equicontinuity } \Rightarrow \hbox{ distality }\Rightarrow \hbox{ entropy 0.}$$
As distality, as well as equicontinuity, excludes Li-Yorke chaos, this inverse scale does not follow the scale of partial properties in a symmetric way. On the matter of disjointness, 
remember that u.p.e. systems are disjoint from minimal zero-entropy systems, and weakly mixing systems are disjoint from minimal distal systems, but there can exist no such disjointness theorem involving  second category scrambled sets, as was already remarked above when dealing with completely scrambled systems. The lowest degree of the second scale is not a very secure one. 

\medskip
With sensitivity as the weakest assumption one can draw an overall scale of chaos:
$$\hbox{uniform positive entropy } \Rightarrow \hbox{ weak mixing }\Rightarrow \hbox{ sensitivity,}$$
then a hypothetical partial scale of chaos in which the last definition is left to the reader:
$$\hbox{positive entropy } \Rightarrow \hbox{ partial weak mixing }\Rightarrow \hbox{ `partial sensitivity',}$$
where the second implication is easily checked provided one chooses a sufficiently weak definition of partial sensitivity, and finally
$$\hbox{equicontinuity } \Rightarrow \hbox{ distality }\Rightarrow \hbox{ entropy 0}$$
on the deterministic side. Here the opposition between the two sides is more satisfactory: equicontinuity is really an opposite of sensitivity, and distality an opposite of weak mixing without contradicting sensitivity. Nevertheless there exists no disjointness theorem between equicontinuity and sensitivity: though sensitive a Sturmian system is never disjoint from the associated (minimal equicontinuous) irrational rotation.

\medskip
Of course many other scales may be considered and other chaotic properties introduced. Ideally one would have to draw a complete diagram including all properties that have to do with chaos or order, all implications and disjointness properties between them. This of course cannot be done in a legible way. The scales I drew partly illustrate what can be done when starting from the two elementary chaotic properties we picked up, Li-Yorke chaos and sensitivity.

\section{Classes of dynamical systems}\label{classes} 

It is striking how chaos may take up different features when one restricts oneself to particular classes of dynamical systems. This may be deceptive. The existence of a dense set of periodic points is implicitly or explicitly considered as a chaotic property in many articles, some of which I am a co-author of. As explained above this cannot be considered true theoretically; and there is a family of easily computable systems,  cellular automata, inside which it is decidedly not true. Here we examine specific implications between chaotic properties and mention counter-examples, all of them occurring in restricted classes. 

\subsection{Interval maps and dimensions 1 and 2}\label{1-2}

It is in the field of interval maps that transition to chaos was first observed among one-parameter families of unimodal transformations. The picture is striking. As the parameter $a$ increases there appear points with period 2, then 4, and so on until suddenly at some value of $a$ the system has infinitely many periodic points and becomes Li-Yorke chaotic; just after this value the transformation becomes sensitive with positive entropy and some degree of transitivity. The behaviour of unimodal transformations is extremely valuable for the understanding of the sudden appearance of disorder in systems. Nevertheless, in abstract terms the conception that all features of chaos arise simultaneously at some threshold is misleading. When considering the whole field of topological dynamics one has the feeling that there ways of progressing less brutally in the scales of chaoticity. 

In the rather large class of interval maps the dynamics varies within narrow limits and various chaotic topological properties are strongly correlated. They have been investigated intensively since the sixties, as can be seen in \cite{Bl}, and the strong relationships existing between chaotic properties are reviewed in \cite{Ru}. Because of the existing literature, and as properties of interval maps are rather widely known among researchers in several fields, I shall not develop this topic here. Let me simply mention that properties that are not always chaotic outside this field, like transitivity or dense periodic orbits, become strong chaoticity criteria inside. For instance transitivity implies Li-Yorke chaos, sensitivity, positive entropy, dense periodic orbits, and is not very far from weak mixing.  Conversely when an interval map has positive entropy there exists  some sub-interval on which  its restriction is transitive. 

Chaos among circle maps is not very different, in spite of the existence of irrational rotations. 
My information is not sufficient about maps on trees or graphs in general, and they are still being investigated. But, though results are much harder to prove in those settings, all those I heard about seem to confirm that the situation is not very different from the one prevailing among interval maps. 
 
\bigskip
There do not seem to be many general results about maps of the square, but the variety of behaviours is certainly huge, unlike what happens in dimension 1.

In many examples of strange chaotic and non-chaotic attractors the underlying space has dimension 2.  Triangular maps of the square, that is, continuous maps of the form $F(x,y)=(f(x), g(x,y))$, $x,\ y \in [0,1]$, have been used to obtain examples with various properties. Often enough the purpose is to show that properties that never occur simultaneously among interval maps may do so in dimension 2. 
 
\subsection{Dimension 0: symbolic dynamics}\label{symb}

Symbolic dynamics is the study of the action of the shift on sets of infinite sequences of symbols. Precisely, let $A$ be a finite alphabet; consider the set $A^\Z$ (sometimes $A^\N$) of all bi-infinite sequences of letters of $A$.
Symbolic dynamical systems or {\it subshifts} are closed shift-invariant subsets of $A^\Z$, endowed with the action of the shift $\sigma$: $(\sigma(x))_i = x_{i+1}$. The standard topological setting is the product of discrete topologies on all coordinates on the symbolic space $A^\Z$: then $A^\Z$ and $A^\N$ are Cantor sets and $\sigma$ is a continuous self-map, invertible in the former case only. 

Putting aside the fact that all transitive subshifts are sensitive, except when they are reduced to one periodic orbit, there are very few general properties of subshifts \cite{LM}. This makes symbolic dynamics an extraordinary mine of examples and counter-examples for topological dynamics, in particular when considering chaos. This has been known for a long time. Here are some examples. \\
1. A Sturmian system is a symbolic representation of the rotation of the unit circle by an irrational number $\alpha$. It is obtained by coding the trajectories of points of the circle according to the position of each iterate with regard to the partition $\{[0,1-\alpha], [\alpha, 0]\}$. It is a one-to-one extension of the rotation, except on one orbit on which it is 2-to-one, which provides a very close joining of the two systems. Sturmian shifts are Auslander-Yorke but not Devaney chaotic, sensitive without Li-Yorke chaos, and they are one further illustration of the fact that transitivity cannot be considered a chaotic property. Though having an overall chaotic property,  sensitivity, they have a non-trivial equicontinuous factor. \\
2.  Toeplitz  systems exhibit a large sample of behaviours with respect to chaos. All of them are extensions of adding machines, which are equicontinuous systems.  But some have positive entropy \cite{Gr,Wi} \\
3. Substitution minimal systems may be bounded-to-one extensions of equicontinuous systems, but others are much more chaotic. The first example of a weakly but not strongly mixing system was a substitutive one, the Chac\'on shift, implicitly introduced in \cite{C}; it is at the same time weakly mixing, therefore both sensitive and Li-Yorke chaotic, and minimal, hence without periodic points. Its entropy is 0.

The Chac\'on shift and many other examples provide further reasons for considering Devaney's definition as a requirement for chaos among restricted classes of dynamical systems only. 

The next two examples show that in general Devaney chaos has no connection with positive entropy.
\\
4. A weakly mixing subshift with dense periodic orbits and entropy 0 was first described in \cite{W}. \\
5. The first explicit examples of minimal systems with positive entropy, introduced in \cite{Gr}, are symbolic systems. 

Other subshifts with properties that may be significant from the point of view of chaos have been described more recently: \\
6. in \cite{HY4} one finds an explicit construction of a system that is not weakly mixing but scattering (an overall property which I preferred not to introduce here; it is weaker than weak mixing but still implies transitivity and Li-Yorke chaos \cite{HY2}); \\
7. in \cite{BH} a Devaney chaotic system which is not partly weakly mixing; \\ 8. finally a strongly mixing system containing a residual scrambled set without being completely scrambled is described in \cite{BHS}; it happens to be a non-minimal substitutive system. 

Example 1 (Section \ref{results}) is not altogether symbolic but derived from a symbolic system. Maybe I am forgetting several other meaningful constructions from symbolic dynamics. 

Because of the remarkable plasticity of the theory such counter-examples have usually a strong abstract character: it is not likely that they will soon become models for physical phenomena. But of course this might happen unexpectedly any day. 

A subshift is entirely described by its associated language, that is, the set of words that occur as blocks of coordinates. All corresponding languages can be described with the help of automata with various degrees of complexity. From the point of view of chaos the landscape changes dramatically when one considers classes like  sofic subshifts, those that are recognised by finite automata, or substitution minimal systems, which are generated  finitely. Subshifts belonging to one of those classes have many common properties, which are radically different in the two cases. Sofic systems, except almost trivial ones, have strong chaotic features and their few possible non-chaotic features are completely described; substitution minimal systems are rather deterministic  and one can often tell whether they have chaotic properties or not.

\subsection{Dimension 0 or $\infty$: cellular automata and their topologies}\label{CA}

Given a finite set $A$ and a positive integer $d$, cellular automata (CA) are  continuous shift-commuting homomorphisms of the symbolic space of dimension $d$, $A^{\Z^d}$, endowed with the semi-group of shifts. For $d>1$ there are many results, but the theory of CA is considerably more developed and the constraints are stronger when $d=1$. For this reason it is the only case we consider here. One-dimensional CA form a countable family of maps, for which computer simulations are easy. They nevertheless exhibit a rather large variety of heuristic behaviours. As yet much less is known about them than about interval maps. 

Most of the specificities of CA acting on $A^\Z$ endowed with the standard topology of symbolic dynamics are listed in \cite[Chapter 5]{Ku}. Let us mention three facts, among others that are true in this setting but not for all topological dynamical systems, even in the compact case: \\  1. transitivity implies sensitivity, which makes it a chaotic property among CA \cite{Ku}; \\ 2. a cellular automaton which is not sensitive has equicontinuity points \cite{Ku}--this divides the family into two classes, one rather deterministic and the other rather chaotic;\\ 3. there are no minimal CA, since any CA has at least one periodic point. 

In the field of CA there does not seem to be any correlation between the complexity of the maximal attractor and the dynamics. One knows several simple cellular automata having a complicated, sometimes non-described, maximal attractor while the dynamics on it, or sometimes on its most significant subset, is trivial. 

Among cellular automata, having a dense set of periodic orbits cannot be considered a chaotic property. In \cite{BK} the authors prove that right- or left-closing CA have dense periodic orbits; those families which I do not want to define here contain the class of positively expansive CA, which are usually considered the most chaotic on several accounts--I leave it to interested readers to find the definitions and consequences of right-closing, left-closing and positive expansivity in \cite{Ku}. But CA with equicontinuity points also have dense periodic orbits, and they are certainly the least chaotic \cite{BT}! Here the property of having dense periodic orbits happens along with chaotic properties but holds systematically for the least chaotic class of CA. Have all surjective CA a dense set of periodic orbits? Today the question is still open.

\medskip
There are many other topological spaces, all of them closely related to $A^\Z$ but very different from it, on which CA can be proved to act. One I shall not mention any more, in spite of its interest, is the Weyl space  \cite{BFK,BCF}. Here I address only two, the standard Cantor space and the Besicovitch space. The latter is defined as follows. For $x,\ y \in A^\Z$ consider the pseudo-distance $$d_B(x,y) = \limsup_{n \to \infty} \frac{1}{2n+1}\ \hbox{Card } \{i \in \Z \ | -n\le i\le n,\ x_i \ne y_i\};$$ it is indeed a pseudo-distance because $d_B(x,x) = 0$ for $x \in A^\Z$ and the triangle inequality holds, but two distinct points are at $d_B$-pseudo-distance 0 if and only if the upper density of the coordinates at which they disagree is 0. For instance when $x$ and $y$ have a finite number of distinct coordinates then $d_B(x,y) = 0$. 

On the quotient space $X_B$ of $A^\Z$ by the equivalence relation $\sim$ defined by $x\sim y$ if and only if $d_B(x,y) =0$, $d_B$ generates a distance, also denoted by $d_B$. The corresponding topology is not compact, and the space $X_B$, called the Besicovitch space, has infinite dimension \cite{BFK}. On $X_B$ the distance $d_B$ is shift-invariant, which of course implies that the shift is a continuous map on this space. Actually all cellular automata act continuously on $X_B$ \cite{CFM}. 

The really surprising fact (not so surprising when considering that $A^\Z$ is compact and $X_B$ is not, which means that the two spaces, however close they look, are deeply different) is that the dynamics of one CA may vary considerably, when one considers it first in the standard topology of $A^\Z$ and then on the Besicovitch space. The shift $\sigma$ is mixing and has positive entropy when acting on $A^\Z$, but it is evidently an isometry on $X_B$; it is a simple exercise to show that an isometry is never weakly mixing and that its entropy is 0, even in a non-compact space. This is not the only CA the properties of which change radically. On the other hand, if a CA has equicontinuity points in $A^\Z$ it has the same property in $X_B$, and if it sensitive in $X_B$ it must be sensitive in $A^\Z$.

On $X_B$  it is hardly possible to talk of chaotic CA. Positively expansive CA, which exist and are strongly chaotic on $A^\Z$, do not exist on $X_B$  \cite{BFK}. There are some sensitive CA on the latter space \cite{BFK,BCF}, but none is transitive \cite{BCF}. Another observation is that  all CA that are already known to be sensitive on $X_B$ do not possess a dense set of periodic points in the same space \cite{BF}. But on the topic of the Besicovitch topology several significant questions have not yet been solved or even addressed. 

\medskip 
These two topological spaces on which CA act can be seen as two different ways of looking at them. In the usual Cantor topology they exhibit a large spectrum of deterministic and chaotic behaviours. In the Besicovitch topology they never look very chaotic; in some instances this non-chaoticity is closely related to the simplicity of computations they generate.  This shows how much chaoticity may depend on the topological context. 

\section{Conclusion and questions}
One of the aims of this paper is to show that one cannot express the essence of chaos in a few lines of mathematics. 

Among the "classical" definitions of chaos that are examined here, Li-Yorke chaos is a feature all systems that are usually considered chaotic possess but it is not a completely satisfactory property, and not enough to ensure chaoticity. The significance of sensitivity is more or less the same, except that it is an overall property. Again because they are overall definitions Auslander-Yorke and Devaney are not mainstream chaotic properties. None of the three first goes very deep into the structure of dynamical systems. Though significant for some types of models  Devaney is slightly out of the way from a purely topological point of view. Other requirements, like weak mixing or positive entropy, go deeper into the structure of the system; they are too strong and abstract to be considered `definitions of chaos'. Intuitively they are less obviously chaotic properties than sensitivity, but there is a long-standing interpretation of entropy as a measure of chaos in physics and information theory. 

Another purpose of the article is to introduce chaotic and order properties, compare them whenever that is possible, find out how opposite a given order property and a given chaotic one are, determine how convenient they are in the theory. This amounts to selecting some notions of topological dynamics and all relevant results in order to make them into a kind of abstract theory of chaos. I hope to have shown that this can be done without completely forgetting the heuristic aim. 
As a consequence it is often possible to say that in this theory one system is more chaotic than another. Of course the theory of chaos that is sketched here, especially in Subsection \ref{degrees}, is tentative, because of the reference to experiment, and incomplete since there are other properties of topological dynamics that ought to be examined. 

Another point  deserves to be noticed. Among chaotic and deterministic notions some are partial and some overall. Unlike several corresponding overall notions partial order and partial chaos may co-exist inside the same dynamical system.  It seems to me that this situation appears in genomic models and may also arise in the formal study of musical works. 

The existence of a gap between the heuristic idea of chaos and the mathematical view cannot be avoided. Devaney's definition is certainly more valuable for physicists than for mathematicians, and may stay so for a long time, perhaps forever. The theoretical point of view is necessary because if one wants to translate an experimental notion into mathematical terms this must be done consistently. 

Well, then how could one make an abstract theory of topological chaos less purely theoretical? This is perhaps the most important challenge for the understanding of chaos and order. An interesting answer might arise from a discussion of sensitivity. Broadly speaking sensitivity means that  there exists some $n$ such that $d(F^n(x), f^n(y)) < \eta$ while $x$ and $y$ are close to each other up to some $\epsilon$. But what are the features of the set of all such values of $n$ given $\epsilon$? How do these features vary when the properties of the sensitive system $(X,T)$ vary? One would like to say that Sturmian systems are not very sensitive. One would also like weak mixing or positive entropy to imply some stronger kind of sensitivity. This cannot be done without introducing new definitions. At present I have not got the least idea of what they might be. Perhaps Kolmogorov complexity could help? Long ago Brudno proved that there exists a relation between this complexity and the measure-theoretic entropy of a system \cite{Bru}. But it is possible that completely new ideas would be required. 

\bigskip
\noindent {\bf Questions.}

Here I gather a few questions that have already been posed above, directly or indirectly.

The first two are standard open questions. Remind that it is known that any extension of a system having scrambled pairs has scrambled pairs \cite{B-M}, but Example 1 shows that an extension of a Li-Yorke chaotic system is not Li-Yorke chaotic in general. 

\smallskip \noindent
1. Is partial weak mixing preserved under extensions? \\ 2. As yet the same question for the three existing  definitions of distributional chaos \cite {S} is also unsolved. 

\smallskip
The next question is much less precise and more ambitious. It concerns notions  that are yet undefined  and which I believe would be useful for a mathematical understanding of chaos and order.

\smallskip \noindent
3. Apart from sensitivity, are there computational properties of dynamical systems which are significant in terms of chaos? Is there a way of showing that the iterated image of a point $x \in X$ is more difficult to compute precisely when $(X,f)$ is weakly mixing or when its entropy is positive than when neither of those  properties holds? 

\bigskip
\noindent {\bf Acknowledgements.}
Wen Huang allowed me to include Example 1, which we found during a discussion, in this article: I want to express my gratitude to him. I want to thank the organisers of the Visegrad conference on dynamical systems (\v Strbsk\'e Pleso, June 2007) at which this example was found and several changes in this article were suggested to me when listening to the talks and discussing with participants. I am also grateful to the referee of JDEA, who in his/her turn proposed Example 1.

\end{document}